\documentclass[12pt]{article}
\usepackage{}
\usepackage{mathrsfs}
\usepackage{mathrsfs}
\usepackage{amsfonts}
\usepackage{amsmath}
\usepackage{amssymb}
\usepackage[all]{xy}
\usepackage{xypic}
\usepackage{graphicx}
\usepackage{CJK}
\usepackage{appendix}
\usepackage{indentfirst}

\allowdisplaybreaks

\begin{document}
\baselineskip=15pt
\title{A note on the strong consistency of M-estimates in linear models\thanks{This work is supported by the National Natural
Science Foundation of China (11171001, 11201001, 11426032), the Natural Science Foundation of Anhui Province (1308085QA03, 1408085QA02), the Science Fund for Distinguished Young Scholars of Anhui Province (1508085J06) and Introduction Projects of Anhui University Academic and Technology Leaders.}}
\author{Xinghui Wang, Shuhe Hu\thanks{Corresponding author. E-mail address:
hushuhe@263.net (S.H. Hu);~wangxinghuial@163.com (X.H. Wang).} \\
{\small{\it Department of Statistics, Anhui University, Hefei 230601, P.R. China}}}
\date{}
\maketitle
\begin{minipage}{155mm}
\indent{\bf Abstract}  We improve a known result on the strong consistency of M-estimates of the regression parameters in a linear model for independent and identically distributed random errors under some mild conditions.
\textbf{MR(2010) Subject Classification} 62F12\\\\
\end{minipage}
\vspace{0.3cm}
\setcounter{equation}{0}
We consider the linear model
        \begin{eqnarray}
        Y_i=x_i'\beta_0+e_i, i=1,2,\cdots,n,\label{1}
        \end{eqnarray}
        where $x_i$ are $p\times1$ known design vectors, $\beta_0$  is a $p\times1$ unknown vector of regression coefficients and $\{e_i\}$ are error variables. The M-estimate $\hat{\beta}_n$ of $\beta_0$ is defined by minimizing
        \begin{eqnarray}
        \sum\limits_{i=1}^n\rho(Y_i-x_i'\beta),\label{2}
        \end{eqnarray}
        where $\rho$ is a convex function. Important examples include Huber's estimate with $\rho(x)=(x^2I(|x|\leq c))/2+(c|x|-c^2/2)I(|x|>c)$, $c>0$, where $I(A)$ is the indicator function of the set $A$, the $\mathscr{L}^q$ regression estimate with $\rho(x)=|x|^q$, $1\leq q\leq 2$, and regression quantiles with $\rho(x)=\rho_\alpha(x)=\alpha x^++(1-\alpha)(-x)^+$, $0<\alpha<1$, where $x^+=\max(x,0)$. In particular, if $q=1$ or $\alpha=1/2$, then the minimizer of (\ref{2}) is called the least absolute deviation estimate. There is a substantial amount of work concerning asymptotic properties of M-estimates.

        We assume that $\rho$ is a non-monotonic convex function on $\mathbb{R}$ with right and left derivatives $\psi_+$ and $\psi_-$. Choose $\psi$ such that $\psi_-(u)\leq \psi(u)\leq \psi_+(u)$ for all $u\in\mathbb{R}$. Write $S_n=\sum_{i=1}^nx_ix_i'$, $d_n=\max_{1\leq i\leq n}x_i'S_n^{-1}x_i$, and assume that $S_{n_0}>0$ for some integer $n_0$ and that $n\geq n_0$.

        Zhao (2002) established the following result.\\\\
        \textbf{Theorem A.} Assume that $\{e_i,i\geq1\}$ is a sequence of independent and identically distributed (i.i.d.) random variables. Suppose there exist positive constants $\Delta$, $C_0$, $C_1$ and $\delta\in(0,1]$ such that the following conditions are satisfied:
        \begin{eqnarray}
        \psi(u+h)-\psi(u)\leq C_0~~\textrm{for}~~ h\in(0,\Delta)~~\textrm{and}~~ u\in \mathbb{R},\label{9}
        \end{eqnarray}
        \begin{eqnarray}
        \textrm{E}\psi(e_1)=0, ~~~|\textrm{E}\psi(e_1+u)|\geq C_1|u|~~\textrm{ for}~~ |u|<\Delta,\label{10}
        \end{eqnarray}
        $$d_n=O(n^{-\delta}).$$
        Assume further that $\textrm{E}|\psi(e_1)|^{1/\delta}<\infty$ if $0<\delta<1$ and $\textrm{E}|\psi(e_1)|^q<\infty$ for some $q>1$ if $\delta=1$. Then $\hat{\beta}_n\rightarrow\beta_0$ a.s as $n\rightarrow\infty$.

We will further discuss the strong consistency of $\hat{\beta}_n$ and obtain the following result.\\\\
        \textbf{Theorem 1. }In model (\ref{1}), assume that $\{e_i,i\geq1\}$ is a sequence of i.i.d. random variables. Suppose that conditions (\ref{9}) and (\ref{10}) are satisfied and
        \begin{eqnarray*}
        d_n=O(n^{-1}).\label{11}
        \end{eqnarray*}
        Then $\hat{\beta}_n\rightarrow\beta_0$ a.s as $n\rightarrow\infty$.\\\\
        \emph{Proof} For the technical proof, see the Appendix.\\
        \emph{Remark.} In Theorem A, $\textrm{E}|\psi(e_1)|^q<\infty$ for some $q>1$ is needed for the case $\delta=1$. Theorem 1 is obtained without the condition $\textrm{E}|\psi(e_1)|^q<\infty$ for some $q>1$ in the case $\delta=1$. Theorem 1 improves the result of Theorem A.

\newpage
\appendix
\setcounter{equation}{0}
\section*{Appendix}
\renewcommand{\theequation}{A.\arabic{equation}}

        To prove the main results of the paper, we need the following lemmas.\\
\textbf{Lemma 1} (c.f. Bennett, 1962). Assume that $\{X_n,n\geq1\}$ is a sequence of independent random variables such that $\textrm{E}X_n=0$ and $|X_n|\leq b$ for all $n\geq1$ and some $b>0$. Denote $B_n^2=\sum_{i=1}^n\textrm{E}X_i^2$. Then for all $\varepsilon>0$,
        \begin{eqnarray*}
        P\Big(\Big|\sum\limits_{i=1}^nX_i\Big|>\varepsilon\Big)\leq 2\exp\Big\{-\frac{\varepsilon^2}{2b\varepsilon+2B_n^2}\Big\}.
        \end{eqnarray*}\\
\textbf{Lemma 2} (c.f. Choi and Sung, 1987). Let $\{X_n,n\geq1\}$ be a sequence of i.i.d. random variables with $\textrm{E}X_1=0$. Assume that $\{a_{ni}, 1\leq i\leq n, n\geq1\}$ is an array of constants. If $\max_{1\leq i\leq n}|a_{ni}|=O(n^{-1})$, then
        \begin{eqnarray*}
        \sum\limits_{i=1}^na_{ni}X_i\rightarrow0~~\textrm{as}~~n\rightarrow\infty.
        \end{eqnarray*}\\
\emph{Proof of Theorem 1.} Let $x_{ni}=S_n^{-1/2}x_i$, $\beta_{n0}=S_n^{1/2}\beta_0$ and $\hat{\beta}_n^*=S_n^{1/2}\hat{\beta}_n$. Then
        \begin{eqnarray}
        \sum\limits_{i=1}^nx_{ni}x_{ni}'=I_p,~~\sum\limits_{i=1}^n||x_{ni}||^2=p,~~d_n=\max\limits_{1\leq i\leq n}||x_{ni}||^2,\label{12}
        \end{eqnarray}
        where $I_p$ is the $p\times p$ identity matrix, and $||\cdot||$ is the Euclidean norm on $\mathbb{R}^p$.

        Model (\ref{1}) can be rewritten as
        \begin{eqnarray}
        Y_i=x_{ni}'\beta_{n0}+e_i,~~i=1,2,\cdots,n\label{13}
        \end{eqnarray}
        and
        \begin{eqnarray}
        \sum\limits_{i=1}^n\rho(Y_i-x_{ni}'\hat{\beta}_n^*)=\min\limits_{\beta}\sum\limits_{i=1}^n\rho(Y_i-x_{ni}'\beta).\label{14}
        \end{eqnarray}

        Without loss of generality, we assume that the true parameter $\beta_0=0$ in model (\ref{1}), i.e., $\beta_{n0}=0$ in (\ref{13}). Denote the unit sphere $U=\{\beta: \beta\in \mathbb{R}^p, ||\beta||=1\}$. Let $\varepsilon>0$ be any given constant. Without loss of generality, it can be assumed that $2C_2\varepsilon<\Delta$. Define
        \begin{eqnarray*}
        D_n(\beta)=\sum\limits_{i=1}^n\{\rho(e_i-x_{ni}'\beta)-\rho(e_i)\},~~ \beta\in \mathbb{R}^p
        \end{eqnarray*}
        and
        \begin{eqnarray*}
       D_n(\varepsilon n^{1/2}\gamma)&=&\sum\limits_{i=1}^n\{\rho(e_i-\varepsilon n^{1/2}x_{ni}'\gamma)-\rho(e_i)\}\nonumber\\
        &=&\sum\limits_{i=1}^n\int_0^{w_{ni}'\gamma}\{\psi(e_i+t)-\psi(e_i)\}dt+\sum\limits_{i=1}^nw_{ni}'\gamma \psi(e_i),~~\gamma\in U,\nonumber\\
        &=:&I_{1n}(\gamma)+I_{2n}(\gamma),
        \end{eqnarray*}
        where $w_{ni}=-\varepsilon n^{1/2}x_{ni}$. Hence
        \begin{eqnarray}
        \inf\limits_{\gamma\in U}D_n(\varepsilon n^{1/2}\gamma)\geq \inf\limits_{\gamma\in U}I_{1n}(\gamma)+\inf\limits_{\gamma\in U}I_{2n}(\gamma)\geq\inf\limits_{\gamma\in U}I_{1n}(\gamma)-\sup\limits_{\gamma\in U}|I_{2n}(\gamma)|.\label{17}
        \end{eqnarray}

        We can divide $U$ into $N$ parts, $U_1$, $U_2$, $\cdots$, $U_N$, such that the diameter of each part is less than $n^{-2}$ and $N\leq(2n^2+1)^p$. Let $T_j$ be the smallest close convex set covering $U_j$. For a fixed $T_j$, there are three cases as follows.

        i) $w_{ni}'\gamma\geq0$ for each $\gamma\in T_j$, then there exists a $\gamma_{ij}\in T_j$ such that $w_{ni}'\gamma_{ij}=\inf\{w_{ni}'\gamma: \gamma\in T_j\}$.

        ii) $w_{ni}'\gamma\leq0$ for each $\gamma\in T_j$, then there exists a $\gamma_{ij}\in T_j$ such that $w_{ni}'\gamma_{ij}=\sup\{w_{ni}'\gamma: \gamma\in T_j\}$.

        iii) $w_{ni}'\gamma>0$ for some $\gamma\in T_j$, and $w_{ni}'\gamma<0$ for some $\gamma\in T_j$, then there exists a $\gamma_{ij}\in T_j$ such that $w_{ni}'\gamma_{ij}=0$.

        Write
        \begin{eqnarray*}
        G(t)=\textrm{E}\psi(e_i+t),~~\Psi_i(t)=\psi(e_i+t)-\psi(e_i)-G(t), ~~t\in \mathbb{R}.
        \end{eqnarray*}
        By the monotonicity of $\psi$,
        \begin{eqnarray}
        \inf\limits_{\gamma\in U}I_{1n}(\gamma)&\geq&\inf\limits_{1\leq j\leq N}\inf\limits_{\gamma\in T_j}I_{1n}(\gamma)\geq\inf\limits_{1\leq j\leq N}\sum\limits_{i=1}^n\int_0^{w_{ni}'\gamma_{ij}}\{\psi(e_i+t)-\psi(e_i)\}dt\nonumber\\
        &\geq&\inf\limits_{1\leq j\leq N}\sum\limits_{i=1}^n\int_0^{w_{ni}'\gamma_{ij}}G(t)dt\Big\{1-
        \frac{|\sum_{i=1}^n\int_0^{w_{ni}'\gamma_{ij}}\Psi_i(t)dt|}{\sum_{i=1}^n\int_0^{w_{ni}'\gamma_{ij}}G(t)dt}\Big\}.\label{19}
        \end{eqnarray}
        Let $\gamma\in U_j$ and $\gamma_{ij}\in T_j$. By (\ref{12}) and the definition of $U_j$ and $T_j$, for large enough $n$,
        \begin{eqnarray*}
        \sum\limits_{i=1}^n(x_{ni}'\gamma_{ij})^2- \sum\limits_{i=1}^n(x_{ni}'\gamma)^2&\leq&|\sum\limits_{i=1}^n(x_{ni}'\gamma_{ij})^2- \sum\limits_{i=1}^n(x_{ni}'\gamma)^2|\nonumber\\
        &=&\Big|\sum\limits_{i=1}^n(\gamma_{ij}-\gamma)'x_{ni}x_{ni}'(\gamma_{ij}+\gamma)\Big|\nonumber\\
        &\leq&\sum\limits_{i=1}^n\parallel\gamma_{ij}-\gamma\parallel\parallel x_{ni}\parallel^2(\parallel\gamma_{ij}-\gamma\parallel+2\parallel\gamma\parallel)\nonumber\\
        &\leq&\sum\limits_{i=1}^nn^{-2}(n^{-2}+2)\parallel x_{ni}\parallel^2\nonumber\\
        &\leq&3pn^{-2}<1/2,
        \end{eqnarray*}
        combining (\ref{12}), we obtain that for $1\leq j\leq N$,
        \begin{eqnarray}
        \sum\limits_{i=1}^n(x_{ni}'\gamma)^2&\geq&\sum\limits_{i=1}^n(x_{ni}'\gamma)^2-1/2=\gamma'\sum\limits_{i=1}^nx_{ni}x_{ni}'\gamma-1/2\nonumber\\
        &=&\parallel\gamma\parallel^2-1/2=1/2.\label{21}
        \end{eqnarray}
        By (\ref{11}) and the selection of $\varepsilon$, for large enough $n$ and for $i=1,2,\cdots,n$ and $j=1,2,\cdots,N$,
        \begin{eqnarray*}
        |w_{ni}'\gamma_{ij}|=|\varepsilon n^{1/2}x_{ni}'\gamma_{ij}|\leq C_2\varepsilon\parallel\gamma_{ij}\parallel\leq C_2\varepsilon(1+n^{-2})<2C_2\varepsilon<\Delta.
        \end{eqnarray*}
        It follows that by (\ref{10}) and (\ref{21}),
        \begin{eqnarray*}
        \inf\limits_{1\leq j\leq N}\sum\limits_{i=1}^n\int_0^{w_{ni}'\gamma_{}ij}G(t)dt&\geq& \inf\limits_{1\leq j\leq N}C_1\sum\limits_{i=1}^n\int_0^{w_{ni}'\gamma_{ij}}tdt\nonumber\\
        &=&\inf\limits_{1\leq j\leq N}\frac{C_1}{2}\sum\limits_{i=1}^n(w_{ni}'\gamma_{ij})^2\nonumber\\
        &\geq&\frac{C_1\varepsilon^2n}{4}.
        \end{eqnarray*}
        For $1\leq j\leq N$,, denote $Y_{ni}^{(j)}=\int_0^{w_{ni}'\gamma_{ij}}\Psi_i(t)dt$. By (\ref{9}) and (\ref{21}), one has that
        \begin{eqnarray*}
        |Y_{ni}^{(j)}|\leq 2C_0|w_{ni}'\gamma_{ij}|<4C_0C_2\varepsilon=:C_3.
        \end{eqnarray*}
        For $1\leq j\leq N$, we also obtain that by
        \begin{eqnarray*}
        \sum\limits_{i=1}^n\textrm{Var}(Y_{ni}^{(j)})&\leq&\sum\limits_{i=1}^n\textrm{E}\{\int_0^{w_{ni}'\gamma_{ij}}(\psi(e_i+t)-\psi(e_i))dt\}^2\nonumber\\
        &\leq&C_0\sum\limits_{i=1}^n|w_{ni}'\gamma_{ij}|\int_0^{w_{ni}'\gamma_{ij}}G(t)dt\nonumber\\
        &\leq&\frac{C_3}{2}\sum\limits_{i=1}^n\int_0^{w_{ni}'\gamma_{ij}}G(t)dt.
        \end{eqnarray*}
        Define event $A_{n}$ by
        \begin{eqnarray*}
        A_{n}:=\Big\{\sup\limits_{1\leq j\leq N}\frac{|\sum_{i=1}^n\int_0^{w_{ni}'\gamma_{ij}}\Psi_i(t)dt|}{\sum_{i=1}^n\int_0^{w_{ni}'\gamma_{ij}}G(t)dt}\geq\frac{1}{2}\Big\}.
        \end{eqnarray*}
        Noting the fact that $N\leq (2n^2+1)^p$ and using Lemma 1, we have that
        \begin{eqnarray}
        P(A_n)&\leq&\sum\limits_{j=1}^NP\Big(\Big|\sum_{i=1}^nY_{ni}^{(j)}\Big|\geq\frac{1}{2}\sum_{i=1}^n\int_0^{w_{ni}'\gamma_{ij}}G(t)dt\Big)\nonumber\\
        &\leq&2\sum\limits_{j=1}^N\exp\Big\{-\frac{\frac{1}{4}(\sum_{i=1}^n\int_0^{w_{ni}'\gamma_{ij}}G(t)dt)^2}
        {2C_3\times\frac{1}{2}\sum_{i=1}^n\int_0^{w_{ni}'\gamma_{ij}}G(t)dt+2C_3\times\frac{1}{2}\sum_{i=1}^n\int_0^{w_{ni}'\gamma_{ij}}G(t)dt}\Big\}\nonumber\\
        &=&2\sum\limits_{j=1}^N2\exp\Big\{-\frac{1}{8C_3}\sum_{i=1}^n\int_0^{w_{ni}'\gamma_{ij}}G(t)dt\Big\}\nonumber\\
        &\leq&2(2n^2+1)^p\exp(-C_4n),\nonumber
        \end{eqnarray}
        thus $\sum_{n=1}^\infty P(A_n)<\infty$. By Borel-Cantelli lemma, it follows that $P(A_n, i.o.)=0$. Hence with probability one for large enough $n$,
        \begin{eqnarray*}
        \sup\limits_{1\leq j\leq N}\frac{\sum_{i=1}^n\int_0^{w_{ni}'\gamma_{ij}}\Psi_i(t)dt}{\sum_{i=1}^n\int_0^{w_{ni}'\gamma_{ij}}G(t)dt}<\frac{1}{2},
        \end{eqnarray*}
        which implies by (\ref{19}) that with probability one for large enough $n$,
        \begin{eqnarray}
        \inf\limits_{\gamma\in U}I_{1n}(\gamma)\geq\frac{C_1\varepsilon^2n}{8}.\label{27}
        \end{eqnarray}
        Denote
          \begin{eqnarray*}
          x_{ni}=
          \left(
          \begin{array}{c}
          x_{ni1}\\
          x_{ni2}\\
          \vdots\\
          x_{nip}
          \end{array}
          \right),
          \gamma=
          \left(
          \begin{array}{c}
          \gamma_1\\
          \gamma_2\\
          \vdots\\
          \gamma_p
          \end{array}
           \right),
          ~~  a_{ni}=
          \left\{
          \begin{array}{cc}
          \frac{x_{nik}}{n^{1/2}}, 1\leq i\leq n\\
          0, i\geq n
          \end{array}
          \right. \textrm{for fixed}~~k=1,2,\cdots,p.
        \end{eqnarray*}
         By (\ref{12}) and $d_n=O(n^{-1})$, for fixed $k=1,2,\cdots,N$, $|x_{nik}|\leq \parallel x_{ni}\parallel\leq d_n^{1/2}\leq \sqrt{C_2}n^{-1/2}$ and $\sum_{i=1}^nx_{nik}^2=1$, thus
         \begin{eqnarray*}
         |a_{ni}|\leq \sqrt{C_2}n^{-1}~~\textrm{for}~~n\geq1.
         \end{eqnarray*}
         Using Lemma 2, it follows that
         \begin{eqnarray*}
         n^{-1/2}\sum\limits_{i=1}^nx_{nik}\psi(e_i)=\sum\limits_{i=1}^na_{ni}\psi(e_i)\rightarrow0~~a.s.
         \end{eqnarray*}
         Hence,
         \begin{eqnarray*}
         n^{-1}\sup\limits_{\gamma\in U}|I_{2n}(\gamma)|&=&n^{-1}\sup\limits_{\gamma\in U}|\sum\limits_{i=1}^nw_{ni}'\gamma\psi(e_i)|=n^{-1}\sup\limits_{\gamma\in U}|\sum\limits_{i=1}^n\varepsilon n^{1/2}x_{ni}'\gamma\psi(e_i)|\nonumber\\
         &=&\varepsilon n^{-1/2}\sup\limits_{\gamma\in U}\Big|\sum\limits_{i=1}^n\sum\limits_{k=1}^px_{nik}\gamma_k\psi(e_i)\Big|\nonumber\\
         &=&\varepsilon n^{-1/2}\sup\limits_{\gamma\in U}\Big|\sum\limits_{k=1}^p\Big(\sum\limits_{i=1}^nx_{nik}\psi(e_i)\Big)\gamma_k\Big|\nonumber\\
         &\leq&\varepsilon n^{-1/2}\sup\limits_{\gamma\in U}\sqrt{\sum\limits_{k=1}^p\Big(\sum\limits_{i=1}^nx_{nik}\psi(e_i)\Big)^2}\sqrt{\sum\limits_{k=1}^p\gamma_k^2}\nonumber\\
         &=&\varepsilon\sqrt{\sum\limits_{k=1}^p\Big(n^{-1/2}\sum\limits_{i=1}^nx_{nik}\psi(e_i)\Big)^2}\rightarrow0~~a.s.,
         \end{eqnarray*}
         which implies that with probability one for large enough $n$,
         \begin{eqnarray}
         \sup\limits_{\gamma\in U}|I_{2n}(\gamma)|\leq \frac{C_2\varepsilon^2n}{16}.\label{31}
         \end{eqnarray}
         Combining (\ref{17}) with (\ref{27}) and (\ref{31}), we have that with probability one for large enough $n$,
         \begin{eqnarray*}
         \inf\limits_{\gamma\in U}D_n(\varepsilon n^{1/2}\gamma)\geq \frac{C_2\varepsilon^2n}{16}.
         \end{eqnarray*}
         By the convexity of $D_n(\cdot)$, $D_n(0)=0$ and the definition (\ref{14}) of $\hat{\beta}_n^*$, it follow that
         \begin{eqnarray*}
         \Big\{\inf\limits_{\gamma\in U}D_n(\varepsilon n^{1/2}\gamma)>0\Big\}\subset\{\parallel\hat{\beta}_n^*\parallel\leq \varepsilon n^{1/2}\}.
         \end{eqnarray*}
         Thus for any given $\varepsilon>0$, with probability one for large enough $n$,
         \begin{eqnarray*}
         \parallel\hat{\beta}_n^*\parallel\leq \varepsilon n^{1/2},
         \end{eqnarray*}
         which implies that
         \begin{eqnarray*}
         n^{-1/2}\parallel\hat{\beta}_n^*\parallel\rightarrow0~~a.s.~~\textrm{as}~~n\rightarrow\infty.
         \end{eqnarray*}
        Denote $\zeta(A)$ for the smallest eigenvalue of a positive definite matrices $A$. Based on the result that if $A$ and $B$ are two positive definite matrices of order $p$, then
         \begin{eqnarray*}
         \textrm{tr}(AB)\geq \mu(A)\zeta(B),
         \end{eqnarray*}
         where $\mu(A)$ is the largest eigenvalue of $A$ (see Chen \& Zhao, 1995). Now we take $M=n_0$ and $n>n_0$ so that $S_n\geq S_M>0$. By (\ref{11}), one has that
         \begin{eqnarray*}
         \zeta(S_M)(\zeta(S_M))^{-1}&\leq& \zeta(S_M)\textrm{tr}(S_n^{-1})\leq \textrm{tr}(S_M^{1/2}S_n^{-1}S_M^{1/2})\nonumber\\
         &=&\textrm{tr}(S_n^{-1}S_M)=\sum\limits_{i=1}^Mx_i'S_n^{-1}x_i\leq Md_n\leq MC_2n^{-1},\nonumber
         \end{eqnarray*}
         and there exists a positive constant $C_5$ such that $1\leq C_5n^{-1/2}\zeta(S_n^{1/2})$. It follows that
         \begin{eqnarray*}
         \parallel\hat{\beta}_n\parallel\leq C_5n^{-1/2}\parallel S_n^{1/2}\hat{\beta}_n\parallel=C_5n^{-1/2}\parallel\hat{\beta}_n^*\parallel\rightarrow0~~a.s.~~\textrm{as}~~n\rightarrow\infty.
         \end{eqnarray*}
         The proof of the theorem is completed.\\

\end{document}